\input amstex.tex
\documentstyle{amsppt}
\magnification=\magstep1
\def\Dt{\partial_t}
\def\({\left(}
\def\){\right)}
\def\eb{\varepsilon}
\def\sppan{\operatorname{span}}
\def\R{\Bbb R}
\NoRunningHeads
\document
\topmatter
\title Counterexamples to the regularity of Mane projections and global attractors.
\endtitle
\author Alp Eden${}^1$, Varga Kalanarov${}^2$ and Sergey Zelik${}^3$ \endauthor
\address{${}^1$ Department of mathematics,  Bogazici University,
\newline\indent Bebek 34342, Istanbul, Turkey}
\endaddress
\email{eden\@boun.edu.tr}\endemail

\address{${}^2$ Department of mathematics, Ko{\c c} University,
\newline\indent Rumelifeneri Yolu, Sariyer 34450\newline\indent
Sariyer, Istanbul, Turkey}
\endaddress
\email{vkalantarov\@ku.edu.tr}\endemail
\address{${}^3$ Department of mathematics,
\newline\indent
University of Surrey Guildford,
\newline\indent  GU2 7XH United Kingdom }
\endaddress
\email{S.Zelik\@surrey.ac.uk}\endemail

\abstract We study the global attractors of abstract semilinear parabolic equations and their projections to finite-dimensional planes. It is well-known that the attractor can be embedded into the finite-dimensional inertial manifold if the so-called spectral gap condition is satisfied.
We show that in the case when the spectral gap condition is violated, it is possible to construct the nonlinearity in such way that the corresponding attractor cannot be embedded into any finite-dimensional Log-Lipschitz manifold and, therefore, does not possess any Mane projections with Log-Lipschitz inverse. In addition, we give an example of finitely smooth nonlinearity such that the attractor has finite Hausdorff but infinite fractal dimension.
\endabstract
\endtopmatter
\head Introduction
\endhead
It is well-known that the long-time behavior of many dissipative systems generated by partial differential equations (PDEs) arising in mathematical physics can be described in terms of the so-called global attractors. Being a compact invariant set of the phase space attracting the images of all bounded sets (as time tends to infinity), the global attractor (if exists) contains all of the non-trivial limit dynamics of the system considered on the one hand and, on the other hand, it is essentially smaller than the initial phase space. In particular, in the case where the underlying PDE is considered in a bounded domain of $\Omega\subset\R^n$, this attractor $\Cal A$ often has the finite Hausdorff and fractal dimension, see \cite{2,4,22,34} and the reference therein.
\par
Thus, despite the initial infinite-dimensionality of the phase space (say, $\Phi=L^2(\Omega)$), the limit dynamics is occurred finite dimensional and is equivalent to the appropriate dynamical system (DS) defined on a compact subset of $\R^N$. Indeed, if the fractal dimension
$$
 \dim_f \Cal A=\kappa<\infty
\tag0.1
$$
then, due to the Mane projection theorem, see \cite{19}, the orthogonal projector $P$ to a generic plane of dimension $N>2\kappa+1$ is one-to-one on the attractor $\Cal A$ and, by this reason, the restriction of   the solution semigroup $S(t):\Phi\to\Phi$, $t\ge0$ associated with the considered dissipative system to the global attractor $\Cal A$ is {\it homeomorphic} to the semigroup $\widetilde S(t)$ acting on $P\Cal A\subset\R^N$:
$$
\tilde S(t):P\Cal A\to P\Cal A,\ \ \tilde S(t)v:=PS(t)P^{-1}v,\ \ v\in P\Cal A
\tag0.2
$$
Thus, at least on the level of {\it topological dynamics}, there is no essential difference between the infinite-dimensional dissipative systems
(with finite-dimensional global attractors) and the classical finite-dimensional DS, see e.g.,  \cite{13} for more details.
\par
However, only the continuity of the reduced DS $\widetilde S(t)$ guaranteed by the Mane projection theorem is {\it not sufficient} for the more or less effective study of that system since the most part of the ideas and methods of the modern DS theory require the DS considered to be smooth, see \cite{14} and references therein. Therefore, the factual {\it smoothness} of the DS \thetag{0.2} becomes crucial for the effectiveness of the finite-dimensional reduction described above. Since the initial infinite-dimensional DS $S(t)$ is usually smooth, the most natural and straightforward way to gain the regularity of the reduced semigroup $\widetilde S(t)$ is to verify the extra regularity of the inverse Mane projector $P^{-1}$ on the reduced attractor $P\Cal A\subset\R^N$ and, by this reason, establishing the additional smoothness of the inverse Mane projector becomes one of the central problem of the theory, see \cite{28} and references therein.
\par
An ideal situation arises when the initial dissipative system possesses the so-called {\it inertial} manifold. That is a finite-dimensional invariant submanifold $\Cal M$ of the phase space diffeomorphic to $\R^N$ which contains the global attractor $\Cal A$ and possesses the exponential tracking property, see \cite{7,9,20,34} and references therein. In that case, the limit dynamics is described by the reduced system of ODEs on the manifold $\Cal M=\R^N$:
$$
y'=\Cal F(y(t)),\ \ y\in\R^N,
\tag0.3
$$
called the {\it inertial form} of the initial DS and the smoothness of the vector field $\Cal F$ is determined by the factual smoothness of the manifold $\Cal M$, see \cite{6} for the details.
\par
However, being a particular case of a center manifold, to exist the inertial manifold requires some  splitting of the phase space to "fast" and "slow" variables and this usually leads to rather restrictive {\it spectral gap} conditions.
\par
Indeed, let us consider the case of an abstract semilinear parabolic equation
$$
\Dt u+Au=F(u), \ \ u\big|_{t=0}=u_0\in H,
\tag0.4
$$
where $A$ is a self-adjoint positive linear operator in a Hilbert space $H$ with compact inverse
 and eigenvalues $0<\lambda_1\le\lambda_2\le\cdots$ and the nonlinear map $F: H\to H$ is smooth and globally Lipschitz on $H$ with the Lipschitz constant $L$ (for simplicity, we will consider here only the case where the operator $F$ does not decrease smoothness and acts from $H$ to itself, although all the results discussed below can be extended with minor changes to the general case where $F$ is subordinated to $A$). Then, the spectral gap condition for the existence of the $N$-dimensional inertial manifold reads
 $$
 \lambda_{N+1}-\lambda_N>2L,
 \tag0.5
 $$
see \cite{10,20,30} and references therein. In particular, even for the reaction-diffusion systems, such a spectral gap exists only in the case of one spatial dimension (see also \cite{17} for some particular cases in higher dimensions) and without the spectral gap the manifold may not exist, see counterexamples in \cite{18,31} as well as Section 1 below). Note also that even in this ideal case the spectral gap is not large enough to guarantee the $C^2$-regularity of the manifold, so in general, one can only prove the $C^{1+\eb}$-regularity of $\Cal M$ for some small positive $\eb$.  \par
Thus, without the spectral gap condition (and without the inertial manifold), one cannot expect
the projected system \thetag{0.2} to be $C^1$-smooth. On the other side, the improved version of the Mane projection theorem, see \cite{3,8,12,23} guarantees that the inverse projector $P^{-1}$ is {\it H\"older} continuous on $P\Cal A$ with some H\"older exponent $\alpha<1$ for generic $N$-dimensional planes with $N>2\kappa+1$. Therefore, the reduced DS \thetag{0.2} is in a fact {\it H\"older} continuous and, under the natural additional smoothness assumptions on the initial dissipative system, one can make the H\"older exponent $\alpha$ arbitrarily close to one by increasing the dimension $N$, see \cite{12,26,28} and references therein. It is also known that the reduced DS $\widetilde S(t)$ can be described by a system of differential equations \thetag{0.3} in $\R^N$ for the proper H\"older continuous vector field $\Cal F$, see \cite{6}. However, this presentation of the reduced DS in the form of ODEs is not very effective since the H\"older regularity of the vector field $\Cal F$
is not sufficient for the uniqueness and one still needs to use the initial dissipative system in order to chose the proper unique trajectory of the reduced ODEs \thetag{0.3}. This problem would be resolved if we would be able to verify some extra regularity (for instance, the Lipschitz continuity)
of the inverse Mane projector (which is not forbidden by the $C^1$-counterexamples to the inertial manifold theory mentioned above) and guarantees the uniqueness for the reduced ODEs \thetag{0.3}.
\par
The problem of establishing the extra regularity of the Mane projections beyond the H\"older class is of a big current interest, see \cite{24,25,27,28} and reference therein. In particular, a number of necessary and sufficient conditions for the {\it Lipschitz} continuity of $P^{-1}$ and the vector field $\Cal F$ in the case of abstract semilinear equation \thetag{0.4} has been established by Romanov, see \cite{32,33}. In particular, it will be true iff the operator $A$ is Lipschitz continuous on the attractor as the map from $H$ to $H$:
$$
\|Au_1-Au_2\|_H\le C\|u_1-u_2\|_H,\ \ \forall u_1,u_2\in \Cal A.
\tag0.6
$$
This theory allowed to verify the existence of Lipschitz Mane projections for the class of 1D semilinear reaction-diffusion equations
$$
u_{t}=u_{xx}-f(u,u_x)
\tag0.7
$$
which do not satisfy the spectral gap condition due to the presence of $u_x$ in the nonlinearity and for which the existence of the inertial manifolds is not known.
\par
Nevertheless, we need to mention that verifying condition \thetag{0.6} for more general equations is a very difficult task and, up to now, no effective ways to check it are developed. On the other hand, based on the standard log-convexity and backward uniqueness estimates for parabolic equations, it is possible to verify the weakened version of \thetag{0.6}:
 $$
\|Au_1-Au_2\|_H\le C\|u_1-u_2\|_H\log\frac C{\|u_1-u_2\|_H},\ \ \forall u_1,u_2\in \Cal A,
\tag0.8
$$
see \cite{1,15} (in a fact, for our case where $F\in C^1(H,H)$, the improved version of \thetag{0.8} with the square root of $\log$ also holds).
By this reason, it looks natural to try to relax
 the "too-restrictive" Lipschitz continuity assumption for the Mane projections till the $\log$-Lipschitz continuity, see \cite{27,29}. Indeed, it is easy to see that
 the following $\log$-Lipschitz continuity of the inverse Mane projector in $H^2:=D(A)$:
 $$
 \|P^{-1}(v_1-v_2)\|_{H^2}\le C\|v_1-v_2\|_{H^2}\log\frac C{\|v_1-v_2\|_{H^2}},\ \ \forall v_1,v_2\in P\Cal A\subset\R^N
 \tag0.9
 $$
(the $H^2$-norms in the right-hand side can be replaced by the $H$-norms since all norms are equivalent in the finite-dimensional space)
 is sufficient for verifying the analogous $\log$-Lipschitz continuity of the vector field $\Cal F$ in \thetag{0.3} which, in turn, implies the uniqueness.
\par
It worth to emphasize however that, despite the fact that the idea with $\log$-Lipschitz continuity looks indeed attractive (and many papers devoted to the study the conditions which guarantee it, mainly for abstract sets, not attractors, see \cite{28} and references therein), we do not know  a single example of a dissipative PDE arising in mathematical physics (or even natural example of a semilinear parabolic equation of the form \thetag{0.4}) where this idea would allow to go beyond the Romanov's theory. On the other hand, again to the best of our knowledge, no counterexamples  to the $\log$-Lipschitz continuity \thetag{0.9} were known before.
\par
The main aim of the present paper is to give such counterexamples for the class of abstract semilinear parabolic equations \thetag{0.4}.  In particular, we show that, at least in the class of \thetag{0.4} one cannot expect the $\log$-Lipschitz continuity of the inverse Mane projectors of the global attractor. Namely, we prove the following theorem which can be considered as the main result of the paper.

\proclaim{Theorem 0.1} Let $L_0:=\sup_{n\in\Bbb N}(\lambda_{n+1}-\lambda_n)<\infty$. Then, for every positive $L$ such that
$$
L>\max\{L_0/2,\lambda_2\},
\tag0.10
$$
there exists a $C^\infty$-smooth nonlinear map $F:H\to H$ with the global Lipschitz constant $L$ such that equation \thetag{0.4} possesses a compact global attractor $\Cal A$ and this global attractor cannot be embedded into any finite-dimensional $\log$-Lipschitz manifold and, in particular, does not possess any $\log$-Lipschitz Mane projections.
\endproclaim
We note that the construction of the non-linearity $F$ is strongly based on the refining of known counterexamples to the Floquet theory for abstract semilinear equations with periodic coefficients, see \cite{5,16,21}. We will use this construction in order to embed the "almost" orthogonal sequences with "bad" properties to the associated global attractor, see Sections 2 and 3 for the details.
\par
The paper is organized as follows. In Section 1 we remind the construction of the Romanov's counterexample in order to show that the $C^1$-inertial manifold may not exist if the spectral gap condition is violated. Although this result is a some kind of folk knowledge, we did not find the precise reference in the literature and, since that is important for what follows, we decide to reproduce it here.
\par
In Section 2, we show that the absence of the spectral gap allows to construct a time periodic linear operator $\Phi(t)$ such that all solutions of the abstract parabolic equation
$$
\Dt v+Av=\Phi(t)v
\tag0.11
$$
decay super-exponentially (as $e^{-t^2}$) as $t\to\infty$. Based on this equation, we construct then the nonlinearity $F$ for \thetag{0.4} in such way that the associated global attractor contains two trajectories $u_1(t)$ and $u_2(t)$ such that
$$
\|u_1(t)-u_2(t)\|_H\le Ce^{-\alpha t^2},\ \ t\ge0.
\tag0.12
$$
Due to the Romanov's criterium (see \cite{33}), this excludes the existence of a Lipschitz inertial manifold as well as any finite-dimensional bi-Lipschitz Mane projections.
\par
Section 3 is devoted to the proof of the main Theorem 0.1. In particular, we introduce here the so-called Log-doubling factor of the attractor  (in the spirit the doubling factor and Bouligand dimension, see \cite{24}) which must be finite if the attractor can be embedded into the finite-dimensional Log-Lipschitz manifold. Then, by the proper small perturbations of the example from the previous section in various orthogonal directions, we embed the "almost-orthogonal" sequence with infinite Log-doubling factor into the global attractor. This excludes the existence of Log-Lipschitz Mane projections and shows that, in general, without the inertial manifold, one can expect only the H\"older continuity of the inverse Mane projections.
\par
In Section 4, we construct two slightly different, but related counterexamples. One of them shows that, in the case of {\it finite} smoothness of the non-linearity $F$, the fractal dimension of the attractor may depend on the choice of the phase space $H^s:=D(A^{s/2})$ and, in particular, the example of an infinite fractal, but finite Hausdorff dimension in some $H^s$ is given. The second, potentially more interesting example demonstrates the global attractor $\Cal A$ embedded into the smooth 2D inertial manifold, but such that the "bad" orthogonal projection $Q_2\Cal A$  to the plane of codimension 2 has infinite Log-doubling factor and, therefore, cannot be embedded into any finite-dimensional Log-Lipschitz manifold. In particular, this give rise of an interesting open question:
 \par
 {\it
 Let the attractor $\Cal A$ do not possess any Log-Lipschitz Mane projections. Is it possible to find an "extension" $\bar{\Cal A}\subset\bar H$ of the attractor $\Cal A$ (or/and the DS \thetag{0.4}) to a larger Hilbert space $\bar H$ such that $\Cal A$ is a projection of $\bar{\Cal A}$ on $H$ and the larger attractor $\bar{\Cal A}$ will possess such Mane projections?}
\par
The positive answer on this question would resolve the problem of the absence of Log-Lipschitz Mane projections.
\par
Finally, in Section 5, we discuss some consequences of the obtained results and the related open problems.
\par
{\bf Acknowledgement:} This work is partially supported by TUBITAK, Project number 107T896. The authors also would like to thank D.Turaev for many interesting discussions.

\head Section 1. Abstract parabolic equations without $C^1$-inertial manifolds.
\endhead
In this section we slightly modify the counterexample of \cite{33} in order to show that the spectral gap condition is, in a sense, necessary for the existence of the inertial manifold in the class of abstract parabolic equations:
$$
\Dt u+Au=F(u), \ \ u\big|_{t=0}=u_0\in H.
\tag1.1
$$
Here, as usual, $H$ is an abstract Hilbert space, $A: D(A)\to H$ is a positive self-adjoint operator in $H$ with compact inverse
and $F$ is a nonlinear operator which satisfies some natural assumptions stated below. Let also $\lambda_1\le\lambda_2\le\lambda_n\le\cdots$ be the eigenvalues of the operator $A$ and $\{e_n\}_{n=1}^\infty$ be the corresponding orthonormal basis of eigenvectors.
\par
The main result of this section is the following theorem.
\proclaim{Theorem 1.1} Let the operator $A$ be such that
$$
L_0:=\sup_{n\in\Bbb N}(\lambda_{n+1}-\lambda_n)<\infty.
\tag1.3
$$
Then, for every $L>\max\{\frac12 L_0,\lambda_1\}$, there exists  a  nonlinear smooth operator $F\in C^\infty(H,H)$ such that
\par
1) $F$ is globally Lipschitz on $H$ with Lipschitz constant $L$:
$$
\|F(u)-F(v)\|_{H}\le L\|u-v\|_H, \ \ \forall u,v\in H,
\tag1.4
$$
\par
2) Equation \thetag{1.1} possesses a compact global attractor $\Cal A$ in $H$,
\par
3) There are no finite-dimensional invariant $C^1$-submanifolds in $H$ containing the global attractor $\Cal A$.
\endproclaim
\demo{Proof} To construct the counterexample, we need the following simple Lemma.
\proclaim{Lemma 1.2} Consider the following 2 dimensional linear system of ODEs:
$$
\cases
\frac d{dt}u_n+\lambda_n u_n=Lu_{n+1}\\
\frac{d}{dt}u_{n+1}+\lambda_{n+1}u_{n+1}=-Lu_n.
\endcases
\tag1.5
$$
Then, if $2L>\lambda_{n+1}-\lambda_n$, the associated characteristic equation does not have any real roots.
\endproclaim
\demo{Proof}
Indeed, the characteristic equation reads
$$
\det\(\matrix -\lambda_n-\lambda &L\\-L&-\lambda_{n+1}-\lambda\endmatrix\)=0, \ \ \lambda^2+(\lambda_n+\lambda_{n+1})\lambda+\lambda_n\lambda_{n+1}+L^2=0
\tag1.6
$$
and the roots are $\lambda=\alpha_n\pm i\omega_n$ with
$$
\alpha_n:=-\frac{\lambda_n+\lambda_{n+1}}2,\ \ \omega_n=\frac12\sqrt{4L^2-(\lambda_{n+1}-\lambda_n)^2}>0
\tag1.7
$$
and the lemma is proved.
\enddemo
We are now ready to construct the desired equation \thetag{1.1}. According to \cite{33}, it is enough to find to  equilibria $u_+$ and $u_-$ of
problem \thetag{1.1} such that the spectrum $\sigma(-A+F'(u_-))$ does not possess any real eigenvalues, but the spectrum $\sigma(-A+F'(u_+))$ possesses {\it exactly}
one real eigenvalue which is {\it unstable} (and all other eigenvalues are complex-conjugate).
\par
Indeed, assume that the $C^1$-invariant manifold $\Cal M$ exists and $\dim \Cal M=n$. Then, from the invariance, we conclude that
$$
\sigma_{\Cal M}(u_\pm):=\sigma\((-A+F(u_\pm))\big|_{\Cal T\Cal M(u_\pm)}\)\subset \sigma(-A+F(u_\pm))
\tag1.8
$$
where $\Cal T\Cal M(u_\pm)$ are  tangent planes to $\Cal M$ at $u=u_\pm$ (which belong to the invariant manifold  since these equilibria belong to the attractor). In particular, since the equation is real-valued, analyzing the equilibrium $u_-$, we see that $\dim M$ must be {\it even}
(if $\lambda\in \sigma_{\Cal M}(u_-)$ then $\bar\lambda\in\sigma_{\Cal M}(u_-)$ and $\lambda\ne\bar\lambda$ since there are no real eigenvalues).
\par
Let us now consider the equilibrium $u_+$ here we have exactly one real eigenvalue which should belong to $\sigma_{\Cal M}(u_+)$ since the unstable manifold of $u_+$ belongs to the attractor. Since all other eigenvalues must be complex-conjugate, we conclude that the dimension of $\Cal M$ must be odd. This contradiction proves that the $C^1$-invariant manifold does not exist.
\par
Let us now fix two equilibria $u_\pm=\pm N e_1$ where $N$ is a sufficiently large number and define the maps $F^\pm(u)$ via the coordinates
$F^\pm_n(u):=(F(u),e_n)$, $u=\sum_{n=1}^\infty u_ne_n$,
in the orthonormal basis $\{e_n\}$ as follows:
$$
\multline
F_1^-(u):=-\lambda_1N+Lu_2,\ F^-_2(u):=-L(u_1+N),\\ F_{2n+1}^-(u):=Lu_{2n+2},\ F_{2n+2}^-(u):=-Lu_{2n+1},\ n\ge1
\endmultline
\tag1.9
$$
and
$$
\multline
F_1^+(u):=\lambda_1N+L(u_1-N),\\ F^+_{2n}(u):=-Lu_{2n+1},\ F_{2n+1}^+(u):=Lu_{2n},\ n\ge1.
\endmultline
\tag1.10
$$
Then, both $F^-$ and $F^+$ are smooth and Lipschitz continuous with Lipschitz constant $L$ and we may construct the smooth nonlinear map $F(u)$ such that
$$
F(u)\equiv F^\pm(u) \ \text{if $u$ is close to $u_\pm$}
\tag1.11
$$
and the global Lipshitz constant of $F$ is less than $L+\eb$ (where $\eb=\eb(N)$ can be made arbitrarily small by increasing $N$, see the next section for the more detailed construction of this map). Finally, we may cut-off the nonlinearity $F$ outside of a large ball making it dissipative without expanding the Lipschitz norm (see \cite{34} for the details). This guarantees the existence of a global attractor $\Cal A$.
\par
Let find now the spectrum $\sigma(-A+F'(u_\pm))$ at equilibria $u_\pm$. Indeed, according to the construction of $F$, the linearization of \thetag{1.1}
at $u=u_-$ looks like
$$
\frac d{dt}v_{2n-1}=-\lambda_{2n-1} v_{2n-1}+Lv_{2n},\ \ \frac d{dt}v_{2n}=-\lambda_{2n}v_{2n}-Lv_{2n-1}, \ \ n=1,2,\cdots
\tag1.12
$$
and, due to condition $L>\frac12 L_0$, there are no real eigenvalues in $\sigma(-A+F'(u_\pm))$, see Lemma 1.1.
\par
In contrast to that, the linearization near $u=u_+$ reads
$$
\multline
\frac d{dt} v_1=(L-\lambda_1)v_1,\\
\frac d{dt}v_{2n}=-\lambda_{2n} v_{2n}+Lv_{2n+1},\ \ \frac d{dt}v_{2n+1}=-\lambda_{2n+1}v_{2n+1}-Lv_{2n}, \ \ n=1,2,\cdots
\endmultline
\tag1.13
$$
Thus, since $L>\lambda_1$, we have exactly one positive eigenvalue and all other eigenvalues are complex conjugate by Lemma 1.2. This shows the absence of any finite-dimensional $C^1$-invariant manifold for that equation and finishes the proof of the theorem.
\enddemo
\remark{Remark 1.1} Recall that the condition $L<\frac12 L_0$ is a classical spectral gap condition
 which guarantees the existence of a finite-dimensional inertial manifold for equation \thetag{1.1}, see \cite{20,30}. The extra condition $L>\lambda_1$ is necessary in order to have the instability in equation \thetag{1.1}, otherwise the nonlinearity will be not strong enough to compensate the dissipativity of $A$ and the attractor will consist of a single exponentially stable point (which can be naturally treated as a zero dimensional inertial manifold for that problem). Thus, the proved theorem shows the sharpness of the spectral gap condition in the class of abstract parabolic equations and $C^1$-manifolds.
\endremark
\comment
\head Example II. Parabolic attractor without bi-Lipschitz Mane projections
\endhead
In this section, we refine the counterexample from the previous section to show that without the spectral gap condition, not only $C^1$, but also
Lipschitz invariant manifolds containing the attractor may not exist. Actually, we will find to trajectories $u$ and $v$ on the attractor such that
$$
\|u(t)-v(t)\|_H\le Ce^{-\kappa t^2}, \ \ t\ge0
\tag14
$$
for some positive $C$ and $\kappa$. This will imply that the attractor cannot be  bi-Lipschitz projected to any finite-dimensional plane and cannot be embedded to any invariant Lipschitz submanifold, see \cite{33} for the details. Thus, the main result of this section is the following theorem.

\proclaim{Theorem 2} Let $A$ be a self-adjoint positive operator with compact inverse acting in a Hilbert space $H$ and let the spectral gap exponent
$L_0<\infty$, see \thetag{3}. Then, for every $L>\frac12 L_0$ such that, in addition, $L>\lambda_2$, there exists  and a smooth nonlinearity $F(u)$ satisfying \thetag{4} such that the corresponding abstract parabolic equation \thetag{1} possesses a global attractor $\Cal A$  which contains at least two trajectories $u(t)$ and $v(t)$ satisfying \thetag{14}.
\endproclaim
\demo{Proof} We first note that, without loss of generality we may assume that
$$
\lambda_{n+1}-\lambda_n\ge c_0
\tag15
$$
for all $n\ge0$ (otherwise, we will just use a subsequence $e_{n_k}$ of all eigenvectors $e_n$ of $A$ setting $F_n(u)\equiv0$ if $n\ne n_k$).
\par
The desired counterexample will use the spectral structure from the previous section. However, in addition, we need to build up the recurrent mechanism which will allow the trajectory to move from the neighborhood of $u_+$ to the neighborhood of $u_-$ and back infinitely many times. To this end, the previous infinite-dimensional system will be coupled with the system of 2 ODEs on a plane $(x,y)$:
$$
\frac d{dt}x=f(x,y),\ \ \frac d{dt} y=g(x,y)
\tag16
$$
which is assumed to satisfy the following properties:
\par
1) There are at least 3 equilibria $U_-=(-N,0)$, $U_0=(0,0)$ and $U_+=(N,0)$ such that $U_\pm$ are the saddles and $U_0$ is an unstable focus.
\par
2) There is a stable periodic orbit $\Gamma$ which spends "most" time near saddles $U_+$ and $U_-$ (this will be true, e.g., if $\Gamma$ is obtained by splitting the homoclinic structure consisting of two heteroclinic orbits from $U_-$ to $U_+$ and back.
\par
3) The unstable point $U_0$ belongs to the interior of the bounded domain restricted by $\Gamma$ and all trajectories staring from this interior tend to $U_0$ as time tends to $-\infty$. Without loss of generality, we may assume that the trajectories in this area rotate anti-clockwise.
\par
4) The nonlinearities $f$ and $g$ are cut off outside of a large ball in order to have the dissipativity.
\par
It is not difficult to see that such a planar system exists (even in the class of analytic nonlinearities $f$ and $g$).
\par
We now couple this planar system with the infinite-dimensional system on $u=\sum_{n=1}^\infty u_ne_n$ similar to the one, considered in the previous section. To this end, we introduce two linear operators $F^+$ and $F^-$ in $H$ as follows
$$
\multline
F^-_{2n-1}u=\frac12(\lambda_{2n-1}-\lambda_{2n})u_{2n-1}+(L+\eb_n^-)u_{2n},\\ F^-_{2n}u=-\frac12(\lambda_{2n-1}-\lambda_{2n})u_{2n}-(L+\eb_n^-)u_{2n-1},\ \ n\in\Bbb N
\endmultline
\tag17
$$
and
$$
\multline
F^+_1u=(\lambda_1+1)u_1,\ F^+_{2n}u=\frac12(\lambda_{2n}-\lambda_{2n+1})u_{2n}+(L+\eb_n^+)u_{2n+1},\\
 F^+_{2n+1}u=-\frac12(\lambda_{2n}-\lambda_{2n+1})u_{2n+1}-(L+\eb_n^+)u_{2n},\ n\in\Bbb N,
\endmultline
\tag18
$$
where $L$ is a fixed number satisfying $L>c_0>0$ and the sequences of {\it small} numbers $\eb_n^\pm$ will be specified below.
\par
Let also $\theta\in C^\infty_0(\Bbb R)$ be a smooth cut-off function  such that $\theta(x)\equiv1$, $x\in[N/2,N]$
 and $\theta(x)\equiv0$ for $x\le N/4$ and let
$$
F(x)u:=\theta(x)F^+u+\theta(-x)F^-u.
\tag19
$$
Note that the equilibrium $(U_+,0)$ has the two dimensional unstable manifold containing in particular, the direction of $e_1$. However, up to the moment the dynamics on the $(x,y)$ plane is not coupled with the infinite-dimensional variable $u$. In the next step, we couple it with the 3rd direction $e_1$ by adding the perturbation $R(x,y,u_1)$ which acts only in the directions of $x$, $y$ and $u_1$ and satisfy the following properties:
\par
1) $R\equiv 0$ if $x<0$ or $|u_1|<\delta$ (for some sufficiently small $\delta$)
\par
2) The unstable trajectory $\xi(t)=(x(t),y(t),u_1(t))$ of the equilibrium $(U_+,0)$ which leave the equilibrium along the direction $e_1$ returns to the small neighborhood of the circle $\Gamma$ at some point $x(0)=0$, $y(0)>0$ and $|u_1(0)|\le\delta$.
\par
3) The projection $(x(0),y(0))$ of the point $\xi(0)$ to the $(x,y)$-plane belongs to the {\it interior} of the domain restricted by the circle $\Gamma$.
\par
Again, it is not difficult to construct such a perturbation $R$ in the class of $C^\infty$-maps.
\par
Note that, since $(x(0),y(0))$ belongs to the interior of $\Gamma$, the trajectory $(\bar x(t),\bar y(t))$ of the planar system such that
$\bar x(0)=x(0)$ and $\bar y(0)=y(0)$ can be extended for all $t\in\Bbb R$ as a bounded trajectory (which tends to $(0,0)$ as $t\to-\infty$ and to $\Gamma$ as $t\to+\infty$. Therefore, the expanded trajectory $(\bar x(t),\bar y(t),0)$ solves
$$
\multline
\frac d{dt} x=f(x,y)+R_1(x,y,u_1),\ \frac d{dt} y(t)=g(x,y)+R_2(x,y,u_1),\\  \frac d{dt} u+Au=F(x)u+R_3(x,y,u_1)e_1
\endmultline
\tag20
$$
and belongs to the attractor $\Cal A$ of this equation (since an attractor always contains all bounded trajectories defined for all $\Bbb R$).
Furthermore, the trajectory of \thetag{20} which starts from $(0,y(0),u_1(0),0)$ also belongs to the attractor (since it can be extended for negative times by the trajectory $(\xi(t),0)$ which stabilizes as $t\to-\infty$ to $(U_+,0)$ and, therefore, is bounded). In addition, as we will see below,
this trajectory will converge to $(\Gamma,0)$ as $t\to\infty$ and will never leave the strip $|u_1|\le\delta$ for positive time. By this reason, for positive times, this trajectory has a form $(\bar x(t),\bar y(t),u(t))$ where $(\bar x(t),\bar y(t),0)$ is a solution of the planar system belonging to the attractor which has been fixed above and $u(t)\in H$ solves
$$
\frac d{dt} u+A u=F(\bar x(t))u,\ \ u\big|_{t=0}=(u_1(0),0),\ \ t\ge0.
\tag21
$$
We now show that the solution of this linear non-autonomous (close to time-periodic one) equation decays super-exponentially as time tends
to infinity under the proper choice of small parameters $\eb_n^\pm$, see \thetag{17} and \thetag{18}. To this end, let $T_0=0<T_1<T_2<\cdots$ be
a sequence of subsequent solutions of the equation $\bar x(t)=0$. Then, since the trajectory $(\bar x(t),\bar y(t))$ is close to the circle $\Gamma$ and $\Gamma$ spends "most" time near $U_\pm$, we conclude that $\tau_k:=T_{k}-T_{k-1}$ satisfy
$$
T\le \tau_k\le 2T
\tag22
$$
for some {\it large} $T>0$. We claim that it is possible to fix small correctors $\eb_n^\pm$ in such way that
$$
u_1(t)\equiv u_2(t)\equiv\cdots\equiv u_n(t)\equiv0 \ \ \text{for $t\ge T_n$}
\tag23
$$
which will give the desired super-exponential decay rate of $u(t)$. Note that, due to the structure of equation \thetag{21}, in order to verify \thetag{22}, it is enough to check that
$$
u_{n}(T_{n-1}) \ \text{ is the only non-zero component of $u(T_{n-1})$, $n=1,2,\cdots$}
\tag24
$$
Indeed, it is true for $n=1$ (and $T_0=0$). On the time interval $t\in[T_0,T_1]$ $\bar x(t)\le0$ and $F(\bar x(t))=\theta(-\bar x(t)F^-u$. By this reason, all components of $u(t)$ except of $u_1$ and $u_2$ will be zero identically and $u_1(t)$ and $u_2(t)$ solve
$$
\multline
\frac d{dt}u_1=-\frac12(\lambda_1+\lambda_2)u_1+\theta(-\bar x(t))(L+\eb_1^-)u_2, \\ \frac d{dt}u_2=-\frac12(\lambda_1+\lambda_2)u_2-\theta(-\bar x(t))(L+\eb_1^-)u_1,
\endmultline
\tag25
$$
for $t\in[T_0,T_1]$. We know that $u_2(0)=0$ and need to fix small $\eb_1^-$ in such way that $u_1(T_1)=0$. To this end, it is natural to write $(u_1,u_2)$ in polar coordinates $u_1=R\cos\varphi$, $u_2=R\sin\varphi$. Then, the equation for $\varphi$ reads
$$
\frac d{dt}\varphi =(L+\eb_1^-)\theta(-\bar x(t)), \ \ \varphi(0)=0,\ \ \varphi(T_1)=\frac\pi2+\pi k.
\tag26
$$
This equation can be solved explicitly which gives
$$
L\int_0^{T_1}\theta(\bar x(t))\,dt+\eb_1^-\int_0^{T_1}\theta(\bar x(t))\,dt=\frac\pi2+\pi k.
\tag27
$$
Since $\tau_1$ is large and $\bar x(t)<-N/2$ where $\theta\equiv1$ for "most" time, we conclude that $\int_0^{T_1}\theta(\bar x(t))\,dt\sim T$ and
equation \thetag{27} can be solved with respect to $\eb_1^-$ and $\eb_1^-\sim\frac\pi T$ is indeed small.
\par
Thus, we have proved that it is possible to fix small $\eb_1^-$ in such way that the only non-zero component of $u(T_1)$ is $u_2(T_1)$. After that, for $t\in[T_1,T_2]$ we are moving to the region where $\bar x(t)>0$ and the variables $u_2$ and $u_3$ are coupled (through the operator $F^+$) analogously
to \thetag{26}. The variable $u_1$ is now decoupled and $u_1(T_1)=0$ implies $u_1(t)\equiv0$ on that interval. Solving the analogue of the boundary value problem \thetag{25} for that case, we find small $\eb_1^+$ such that $u_2(T_2)=0$.
\par
After that, on the interval $t\in[T_2,T_3]$, the variables $u_3$ and $u_4$ are coupled by the analogue of \thetag{25} and they are decoupled from
$u_1$ and $u_2$, so $u_1(T_2)=u_2(T_2)=0$ implies that they vanish identically for $t\in[T_2,T_3]$ and solving the boundary value problem, we find $\eb_2^-$ such that $u_3(T_3)=0$.
\par
Repeating these arguments, we construct the sequences $\eb_n^\pm\sim \frac\pi T$ in such way that  \thetag{23} is satisfied.
\par
We claim that the constructed trajectory $u(t)$ satisfies
$$
\|u(t)\|_{H}\le Ce^{-\kappa t^2}
\tag28
$$
for some positive $C$ and $\kappa$. Indeed, since the decay rate of the $n$-th component $u_n$ is governed by the $n$-th eigenvalue $\lambda_{n}$ and
$\lambda_n\ge c_0 n$ (due to \thetag{15}), we have
$$
\|u(T_{n+1})\|_{H}\le Ce^{-c_0 (n+1) T}\|u(T_n)\|_H
\tag29
$$
\enddemo
\endcomment
\head Section 2. Absence of a Lipschitz inertial manifold.
\endhead
In this section, we refine the counterexample from the previous section to show that without the spectral gap condition, not only $C^1$, but also
Lipschitz invariant manifolds containing the attractor may not exist. Actually, we will find to trajectories $u$ and $v$ on the attractor such that
$$
\|u(t)-v(t)\|_H\le Ce^{-\kappa t^2}, \ \ t\ge0
\tag2.1
$$
for some positive $C$ and $\kappa$. This will imply that the attractor cannot be  bi-Lipschitz projected to any finite-dimensional plane and cannot be embedded to any invariant Lipschitz submanifold, see \cite{33} for the details. Thus, the main result of this section is the following theorem.

\proclaim{Theorem 2.1} Let $A$ be a self-adjoint positive operator with compact inverse acting in a Hilbert space $H$ and let the spectral gap exponent
$L_0<\infty$, see \thetag{1.3}. Then, for every $L>\frac12 L_0$ such that, in addition, $L>\lambda_2$, there exists   a smooth nonlinearity $F(u)$ satisfying \thetag{1.4} such that the corresponding abstract parabolic equation \thetag{1.1} possesses a global attractor $\Cal A$  which contains at least two trajectories $u(t)$ and $v(t)$ satisfying \thetag{2.1}.
\endproclaim
\demo{Proof}
The construction of the nonlinearity $F(u)$ with such properties is based on some modification of the counterexample to the Floquet theory, see \cite{5,16}. Namely, we first construct a time periodic operator $\Phi(t)$ such that
the norm of $\Phi(t)$ is arbitrary close to the spectral gap constant $L_0$ and
all solutions of the equation
$$
\Dt w+Aw=\Phi(t)w
\tag2.2
$$
 decay faster than exponential as time tends to infinity.  Let the scalar time-periodic function $x(t)$ with a  period $2\tau$ be given. We assume that this function is odd: $x(-t)\equiv-x(t)$, such that
 and $x(\tau/2)=-x(-\tau/2):=-N$ are the minimal and maximal values of that function respectively and $x(\tau/2-t)=x(t)$ for all $t$. Mention also that, without loss of generality,
  we may assume that
 $$
 c_2k\le \lambda_k\le c_1 k
 \tag2.3
 $$
 for some positive $c_1$ and $c_2$. Indeed, the absence of the spectral gap ($L_0<\infty$) gives the upper bound for $\lambda_k$ and the lower bound can be achieved just by dropping out the unnecessary modes.
 \par
 The following proposition is crucial for what follows.
 \proclaim{Proposition 2.1} Let the assumptions of Theorem 1.1 hold. Then, for every $L>L_0$ and every periodic function $x(t)$,
 there exists a constant $K>0$ and a smooth map $\Cal R:\, \R\to \Cal L(H,H)$ such that
 $$
 \|\Cal R(x)\|_{\Cal L(H,H)}\le L.
 \tag2.4
 $$
 Moreover, if $\Phi(t):=\Cal R(x(t/K))$ be the $2T$-periodic map (with $T=K\tau$), then the Poincare map $P:=U(2T,0)$ (where $U(t,s)$ is a solution operator of \thetag{2.2}: $u(t)=U(t,s)u(s)$) satisfies the following properties:
 $$
 Pe_{2n-1}=\mu_{2n-1} e_{2n+1},\ n\in\Bbb N,\ \ Pe_{2n}=\mu_{2n-2} e_{2n-2},\ n>1,\ \ Pe_2=\mu_0e_1,
 \tag2.5
 $$
 where the positive multipliers $\mu_0$, $\mu_n$ defined via
 $$
 \mu_n:=e^{-T(\lambda_{n}+2\lambda_{n+1}+\lambda_{n+2})/2},\ \  \mu_0:=e^{-T(2\lambda_{1}+\lambda_{2})/2}.
 \tag2.6
 $$
 In particular, all solutions of \thetag{2.2} decay super-exponentially, namely,
 $$
 \|w(t)\|_{H}\le Ce^{-\beta t^2}\|w(0)\|_H,
 \tag2.7
 $$
 where positive $C$ and $\beta$ are independent of $w(0)\in H$.
 \endproclaim
 \demo{Proof}  In order to simplify the notations, we will denote by $x(t)$ the $2T$-periodic function obtained from the initial one by the time scaling.
 \par
 Introduce a pair of smooth nonnegative cut-off functions $\theta_1(x)$ and $\theta_2(x)$ such that
\par
1) $\theta_1(x)\equiv1$, $x\in[N/2,N]$
 and $\theta_1(x)\equiv0$ for $x\le N/4$;
 \par
2) $\theta_2(x)\equiv 1$ for $x\ge N/4$ and $\theta_2(x)\equiv 0$ for $x\le0$.
\par
Let us define  the linear operators $F_\pm$ as follows:
$$
\multline
F^-_{2n-1}(x)u=\frac12(\lambda_{2n-1}-\lambda_{2n})u_{2n-1}\theta_2(-x)+\eb\theta_1(-x) u_{2n},\\ F^-_{2n}(x)u=-\frac12(\lambda_{2n-1}-\lambda_{2n})u_{2n}\theta_2(-x)-\eb u_{2n-1}\theta_1(-x),\ \ n\in\Bbb N
\endmultline
\tag2.8
$$
and
$$
\multline
F^+_1(x)u=0,\ F^+_{2n}(x)u=\frac12(\lambda_{2n}-\lambda_{2n+1})u_{2n}\theta_2(x)+ \eb u_{2n+1}\theta_1(x),\\
 F^+_{2n+1}(x)u=-\frac12(\lambda_{2n}-\lambda_{2n+1})u_{2n+1}\theta_2(x)- \eb u_{2n}\theta_1(x),\ n\in\Bbb N,
\endmultline
\tag2.9
$$
where $\eb>0$ is a positive number which will be specified later. Finally, we introduce
the desired operator $\Phi(t)$ via
$$
\Phi(x(t))u:=F^+(x(t))u+F^-(x(t))u.
\tag2.10
$$
We claim that there exists a small $\eb=\eb(T)>0$ such that the operator thus defined satisfied all assumptions of the proposition. Indeed, let $P_-$ and $P_+$ be the solution operators which map $w(0)$ to $w(T)$ and $w(T)$ to $w(2T)$ respectively:
$$
P_-:=U(T,0), \ \ P_+:=U(2T,T)
$$
where $U(t,s))$ is the solution operator for \thetag{2.2} from time $s$ to time $t$ ( $u(t):=U(t,s)u(s)$).
 Then, the spaces
$$
V_{n}^-:=\sppan\{e_{2n-1},e_{2n}\} \ \ \text{and}\ \ \ V_{n}^+:=\sppan\{e_{2n},e_{2n+1}\},\ n\in\Bbb N
\tag2.11
$$
are invariant subspaces for the linear maps $P_-$ and $P_+$ respectively.
\par
We need to look at $P_\pm e_n$. To this end, we introduce $T_0$ such that $x(T_0)=N/4$ and note that, by the construction of the cut-off functions,
all $e_n$'s are invariant with respect to the solution maps $U(T_0,0)$, $U(T,T-T_0)$, $U(T+T_0,T)$ and $U(2T,2T-T_0)$, so we only need to study the maps
$U(T-T_0,T_0)$ and $U(2T-T_0,T+T_0)$. On these time intervals the cut-off function $\theta_2(x(t))\equiv1$ and the equations read:
$$
\multline
\frac d{dt}u_{2n-1}=-\frac12(\lambda_{2n-1}+\lambda_{2n})u_{2n-1}+\eb\theta_1(- x(t)) u_{2n}, \\ \frac d{dt}u_{2n}=-\frac12(\lambda_{2n-1}+\lambda_{2n})u_{2n}-\eb\theta_1(- x(t)) u_{2n-1},
\endmultline
\tag2.12
$$
for $t\in[T_0,T-T_0]$ and
$$
\multline
\frac d{dt}u_{2n}=-\frac12(\lambda_{2n}+\lambda_{2n+1})u_{2n}+\eb\theta_1(x(t)) u_{2n+1}, \\ \frac d{dt}u_{2n+1}=-\frac12(\lambda_{2n}+\lambda_{2n+1})u_{2n+1}-\eb\theta_1(x(t)) u_{2n},
\endmultline
\tag2.13
$$
for $t\in[T+T_0,2T-T_0]$. To study these equations, we introduce the polar coordinates
$$
 u_{2n-1}+iu_{2n}=R_n^-e^{i\varphi_n^-},\ \ u_{2n}+iu_{2n+1}=R^+_ne^{i\varphi_n^+}
 \tag2.14
 $$
 for problems \thetag{2.12} and \thetag{2.13} respectively. Then the phases $\varphi_n^\pm$ solve the equations
 $$
\frac d{dt}\varphi_n^- =-\eb\theta_1(-x(t))\ \text{and}\ \ \frac d{dt}\varphi_n^+ =-\eb\theta_1(x(t))
\tag2.15
$$
for $t\in[T_0,T-T_0]$ and $t\in[T+T_0,2T-T_0]$ respectively.
\par
Finally, if we fix
$$
\eb:=-\frac{\pi}{2\int_{T_0}^{T-T_0}\theta_1(-x(t))\,dt}=-\frac{\pi}{2\int_{T+T_0}^{2T-T_0}\theta_1(x(t))\,dt}
\tag2.16
$$
both $U(T-T_0,T_0)$ and $U(2T-T_0,T+T_0)$ restricted to $V_{n}^-$ and $V_{n}^+$ respectively will be compositions of the rotation
 on the angle $\pi/2$ and the proper scaling. Namely,
 $$
 \multline
 U(T-T_0,T_0)e_{2n-1}=e^{-(T-2T_0)/2(\lambda_{2n-1}+\lambda_{2n})}e_{2n},\\
 U(T-T_0,T_0)e_{2n}=e^{-(T-2T_0)/2(\lambda_{2n-1}+\lambda_{2n})}e_{2n-1},\\
 U(2T-T_0,T+T_0)e_{2n}=e^{-(T-2T_0)/2(\lambda_{2n}+\lambda_{2n+1})}e_{2n+1},\\
  U(2T-T_0,T+T_0)e_{2n+1}=e^{-(T-2T_0)/2(\lambda_{2n}+\lambda_{2n+1})}e_{2n}
 \endmultline
 $$
Furthermore, on the time intervals $t\in[0,T_0]$ and $t\in[T-T_0,T_0]$, we have the decoupled equations
$$
\multline
\frac{d}{dt}u_{2n-1}+\lambda_{2n-1}u_{2n-1}=1/2(\lambda_{2n-1}-\lambda_{2n})\theta_2(-x(t))u_{2n-1},\\
\frac{d}{dt}u_{2n}+\lambda_{2n}u_{2n}=-1/2(\lambda_{2n-1}-\lambda_{2n})\theta_2(-x(t))u_{2n}
\endmultline
$$
and, therefore,
$$
\multline
U(T_0,0)e_{2n-1}=e^{-\lambda_{2n-1}T_0}e^{1/2(\lambda_{2n-1}-\lambda_{2n})\int_0^{T_0}\theta_2(-x(t))\,dt}e_{2n-1}\\
U(T_0,0)e_{2n}=e^{-\lambda_{2n}T_0}e^{-1/2(\lambda_{2n-1}-\lambda_{2n})\int_0^{T_0}\theta_2(-x(t))\,dt}e_{2n}\\
U(T,T-T_0)e_{2n-1}=e^{-\lambda_{2n-1}T_0}e^{1/2(\lambda_{2n-1}-\lambda_{2n})\int_0^{T_0}\theta_2(-x(t))\,dt}e_{2n-1}\\
U(T,T-T_0)e_{2n}=e^{-\lambda_{2n}T_0}e^{-1/2(\lambda_{2n-1}-\lambda_{2n})\int_0^{T_0}\theta_2(-x(t))\,dt}e_{2n}.
\endmultline
$$
  Thus, for the operator $P_-=U(T,T-T_0)U(T-T_0,T_0)U(T_0,0)$, we have
 $$
\multline
 P_-e_{2n-1}=e^{-\lambda_{2n-1}T_0}e^{1/2(\lambda_{2n-1}-\lambda_{2n})\int_0^{T_0}\theta_2(-x(t))\,dt} e^{-(T-2T_0)/2(\lambda_{2n-1}+\lambda_{2n})}\times\\\times
 e^{-\lambda_{2n}T_0}e^{-1/2(\lambda_{2n-1}-\lambda_{2n})\int_0^{T_0}\theta_2(-x(t))\,dt}e_{2n}=e^{-T(\lambda_{2n-1}+\lambda_{2n})/2}e_{2n},\\ P_-e_{2n}=e^{-\lambda_{2n}T_0}e^{-1/2(\lambda_{2n-1}-\lambda_{2n})\int_0^{T_0}\theta_2(-x(t))\,dt}
 e^{-(T-2T_0)/2(\lambda_{2n-1}+\lambda_{2n})}\times\\\times
 e^{-\lambda_{2n-1}T_0}e^{1/2(\lambda_{2n-1}-\lambda_{2n})\int_0^{T_0}\theta_2(-x(t))\,dt}
 e_{2n-1}=e^{-T(\lambda_{2n-1}+\lambda_{2n})/2}e_{2n-1}
 \endmultline
 \tag2.17
 $$
 and, analogously, for $P_+=U(2T,2T-T_0)U(2T-T_0,T+T_0)U(T+T_0,T)$,
$$
\multline
 P_+e_1=e^{-T\lambda_1/2} e_1,\  P_+e_{2n}=e^{-T(\lambda_{2n}+\lambda_{2n+1})/2}e_{2n+1},\\ P_+e_{2n+1}=e^{-T(\lambda_{2n}+\lambda_{2n+1})/2}e_{2n},\ \ n\in\Bbb N.
\endmultline
 \tag2.18
 $$
This proves the desired spectral properties of the Poincare map $P:=P_+\circ P_-$ (see \thetag{2.5}). Note that, due to \thetag{2.16}, $\eb\to0$ if $T\to\infty$ and the norm of the operator $\Phi(x(t))$ with $\eb=0$ clearly does not exceed $L_0$. Thus, \thetag{2.4} is proved.
\par
Let us check \thetag{2.7}. To this end, it is sufficient to verify that
$$
\|P^Ne_{2n}\|_H\le Ce^{-\beta_1 N^2}
\tag2.19
$$
uniformly with respect to $n\in\Bbb N$. Indeed, according to \thetag{2.5} and \thetag{2.3},  we have
where the multipliers $\mu^\pm_{n}$ satisfy
$$
 e^{-C_2Tn}\le \mu_n^\pm\le e^{-C_1Tn}.
\tag2.20
$$
and, for $N\ge n$,
$$
\|P^Ne_{2n}\|\le e^{-CT(\sum_{k=0}^n 2k+\sum_{k=0}^{N-n} 2k)}=e^{-CT(n(n+1)+(N-n)(N-n+1))}\le e^{-CT/2 N^2}
$$
since for $N\le n$, \thetag{2.20} is obvious, \thetag{2.7} is verified and Proposition 2.1 is proved.
\enddemo
It is now not difficult to construct the desired counterexample. To this end, we generate the $2\tau$-periodic trajectory $x(t)$ as a solution
of the 2D system of ODEs:
$$
\frac d{dt}x=f(x,y),\ \ \frac d{dt} y=g(x,y)
\tag2.21
$$
with smooth functions $g$ and $f$ (cut off for large $x$ and $y$ in order to have the dissipativity). Then, we consider the coupled system
for $u=(x,y,w)$:
$$
\frac d{dt}x=f(x,y),\ \ \frac d{dt} y=g(x,y),\ \ \Dt w+Aw=\Cal R(x)w
\tag2.22
$$
Obviously, system \thetag{2.22} is of the form \thetag{1.1} (we only need to reserve the first two modes $e_1$ and $e_2$ for $x$ and $y$ and
re-denote $Q_3A$ by $A$). It is also not difficult to see that the Lipschitz norm of the nonlinearity in \thetag{2.22} can be made arbitrarily close to $L>L_0$, but in order to make produce the periodic trajectory $x(t)$, $y(t)$ we should be able to destabilize the first two modes which leads to the extra condition $L>\lambda_2$.
\par
Finally, in order to finish the construction, we need to guarantee that at least one of the trajectories of the form $v(t):=(x(t),y(t),w(t))$, $t\ge0$
with non-zero $w(t)$ belongs to the attractor. To this end, we fix the trajectory $v(t)$ of \thetag{2.22} such that $w(0)=e_1$ and
 $w_1(t):=(w(t),e_1)<1$ for all $t\ge0$. After that, we introduce a smooth coupling $R(x,y,w_1)=(R_1,R_2,R_3)$ and the perturbed version of \thetag{2.22}
$$
\multline
\frac d{dt} x=f(x,y)+R_1(x,y,w_1),\ \frac d{dt} y(t)=g(x,y)+R_2(x,y,w_1),\\  \frac d{dt} w+Aw=\Cal R(x)w+R_3(x,y,w_1)e_1
\endmultline
\tag2.23
$$
such that $R\equiv0$ if $w_1\le2$ and such that, in addition, the 3D system
$$
\multline
\frac d{dt} x=f(x,y)+R_1(x,y,w_1),\ \frac d{dt} y(t)=g(x,y)+R_2(x,y,w_1),\\  \frac d{dt} w_1+\lambda_1w_1=R_3(x,y,w_1)
\endmultline
\tag2.24
$$
possesses a saddle point and $(x(0),y(0),1)$ belongs to its unstable manifold and
 the corresponding trajectory $v(t)=(x(t),y(t),w(t))$ satisfies $x(t)<0$ for $t<0$.
Such smooth coupling obviously exists and, since any unstable manifold belongs to the attractor, the trajectory $v(t)$ now belongs to the attractor. By the same reason, the trajectory $u(t):=(x(t),y(t),0)$ also belongs to the attractor. But, according to Proposition 2.1,
$$
\|u(t)-v(t)\|=\|w(t)\|\le Ce^{-\beta t^2}
\tag2.25
$$
and, according to \cite{33}, the Lipschitz inertial manifold containing the global attractor cannot exist.
It is also not difficult to check that, increasing the period of the orbit $x(t)$ by scaling if necessary, it is possible to make the Lipschitz norm of the whole nonlinearity arbitrarily close to $1/2L_0$.
Theorem 2.1 is proved.
\enddemo

\head  Section 3. Absence of the log-Lipschitz embeddings
\endhead
In this section, we show that under the assumptions of Theorem 2.1, it is possible to constrcuct the smooth nonlinearity $F(u)$ in such way that the corresponding attractor $\Cal A$ will be not embedded to any finite-dimensional log-Lipscitz manifold. We recall that a map $T:X\to Y$ between two metric spaces $X$ and $Y$ is log-Lipschitz ($\gamma$-log-Lipschitz) if there exists $\gamma\in(0,\infty)$ and $C>0$ such that
$$
d(Tx_1,Tx_2)\le C d(x_1,x_2)\(\log\frac C{d(x_1,x_2)}\)^\gamma
\tag3.1
$$
for all $x_1,x_2\in X$. The map $T$ is bi-Log-Lipschitz homeomorphism if both $T$ and $T^{-1}$ are log-lipschitz and $\Cal M$ is a log-Lipschiz manifold over $\Bbb R^N$ if it is locally homeomorphic to $\R^N$ and all coordinate maps are bi-Log-Lipschitz (with the same constants $C$ and $\gamma$).
\par
In order to check the absence of such embeddings, we need to introduce some technical tools.
\proclaim{Definition 3.1} Let $X$ be a compact metric space. Then, every its subset $B\subset X$ is pre-compact and, by the Hausdorff criterium, can be covered by the finite number of $\eb$-balls, for every $\eb>0$. Denote by $N_\eb(B,X)$ the minimal number of $\eb$-balls covering $X$. Then, the fractal (box-counting dimension) of $B$ is defined via
$$
\dim_f(B,X)=\limsup_{\eb\to0}\frac{\log N_\eb(B,X)}{\log\frac1\eb},
\tag3.2
$$
see \cite{28} for the details. Define also the doubling factor $D_\eb(X)$ via
$$
D_\eb(X):=sup_{x\in X}N_{\eb/2}(B(\eb,x))
$$
where $B(r,x)$ is the $r$-ball of $X$ centered at $x\in X$.
\endproclaim
It is not difficult to see that $D_\eb(\eb)\le D_N<\infty$ for $\eb\to0$ if $X$ is a bounded subset of $\R^N$ and, therefore, the condition
$$
D(X):=\sup_{\eb>0}D_\eb(X)<\infty
\tag3.3
$$
is {\it necessary} (but not sufficient) for the possibility to embed $X$ into a finite-di\-men\-sio\-nal {\it Lipschitz} manifold, see e.g., \cite{24,28} for the details. The next lemma gives the analogous necessary condition for the case of Log-Lipschitz manifolds.
\proclaim{Lemma 3.2} Let $X$ be a compact set embedded into the Log-Lipschitz manifold $M$. Then the quantity
$$
\dim_{Log-D}(X):=\limsup_{\eb\to0}\frac{\log D_\eb(X)}{\log \log\frac1\eb}
\tag3.4
$$
is finite:
$$
\dim_{Log-D}(X)<\infty.
\tag3.5
$$
\endproclaim
\demo{Proof} Obviously, we need to check \thetag{3.5} for small $\eb$ only, so we need to estimate the number of $\eb/2$-balls covering the ball $B(\eb,x)$, $x\in X$ where $\eb\ll1$.
 Without loss of generality, we may assume that  $B(\eb,x)$ and all covering $\eb/2$-balls belong to the same coordinate chart.  Then there exists a bijective map $T$ of $V\supset B(\eb,x)$ to an open subset  of $\Bbb R^N$ such that $T$ and $T^{-1}$ satisfy \thetag{3.1}. Thus, according to \thetag{3.1},
$$
T^{-1}B(\eb,x)\subset B(\eb_1,T^{-1}x),\ \ T B(\eb_2,y)\subset B(\eb/2,Ty),\ \ y\in T^{-1}V
\tag3.6
$$
with
$$
\eb_1:=C\eb\(\log \frac C\eb\)^\gamma,\ \ C\eb_2\(\log \frac C{\eb_2}\)^\gamma=\eb/2
\tag3.7
$$
Thus, any covering of the ball $B(\eb_1,T^{-1}x)$ by the $\eb_2$-balls in $\Bbb R^N$ will generate the $\eb/2$ covering of $B(\eb, x)$ by the $\eb/2$-balls  balls $X$. By this reason,
$$
N_{\eb/2}(B(\eb,x))\le N_{\eb_2}(B(\eb_1,T^{-1}x))\le C_1\(\frac{\eb_1}{\eb_2}\)^N,
\tag3.8
$$
where the constant $C_1$ is independent of $\eb_i$ and $T^{-1}x$
due to the scaling and shift invariance of $\Bbb R^N$. From \thetag{3.8} and \thetag{3.7}, we now conclude that
$$
N_{\eb/2}(B(\eb,x))\le C_2\(\log\frac C\eb\)^{2\gamma N}
$$
which implies \thetag{3.5} and finishes the proof of the lemma.
\enddemo
We are now ready to state the main result of the section.
\proclaim{Theorem 3.1} Let the assumptions of Theorem 2.1 hold. Then, for
any $L>\max\{1/2 L_0, \lambda_2\}$ there exists a smooth ($C^\infty$-smooth) nonlinearity $F(u)$ with the global Lipschitz constant $L$ such that the associated equation \thetag{1.1} possesses a compact global attractor $\Cal A$ with infinite Log-Doubling factor:
$$
\dim_{Log-D}(\Cal A)=\infty.
\tag3.9
$$
In particular, $\Cal A$ does not possess any finite-dimensional Log-Lipschitz Mane projections and cannot be embedded into any finite dimensional Log-Lipschitz manifold.
\endproclaim
\demo{Proof} We will base on the construction given in the proof of Theorem 2.1. But now we need not only point $v(0)=(x(0),y(0),e_1)$, but the whole
1D {\it segment} $v_s(0):=(x(0),y(0),se_1)$, $s\in[1,1-\kappa]$, $\kappa<1$, to be in the attractor. Then, we will have the whole family of trajectories $v_s(t)$ on the attractor approaching the trajectory $u(t)=(x(t),y(t),0)$ with the super-exponential speed. To achieve this, we only need to modify the 3D system \thetag{2.24} in such way that there is a saddle with 2D unstable manifold $\Cal M_+$ containing the whole segment $v_s(0)$, $s\in[1-\kappa,1]$. We rest the explicit construction of $R$ as an elementary exercise for the reader and will assume from now on that the system \thetag{2.23} is chosen in a such way that $v_s(t)$ belongs to the attractor $\Cal A$ for any $s\in[1-\kappa,1]$.
\par
As in most part of the counterexamples in the dimension and Mane projections theory, see \cite{24,28} and references therein, our counterexample is based on constructing something close to the orthogonal sequence with "bad" properties. However, in contrast to the cited literature, we now need to find how to embed such sequences into the attractor.
\par
Roughly speaking, we will form the prototype of such  sequence by making the small kicks of $v_s(t)$ near $t=0$ in the orthogonal directions of $e_{2n}$ depending on $s$. Of course, it will be still not the proper sequence since all points on it will have large 1st component $w_1$. Then, we will crucially use the fact that the first component $w_1(t)$ of the solution decays "much faster" than the others, so after the properly chosen time
of evolution, we will see indeed the  "almost" orthogonal sequence with "bad" properties. To be more precise, the following lemma holds.
\proclaim{Lemma 3.3} Let $P$ be the Poincare map associated with problem \thetag{2.2}. Then, for every
 $n\in\Bbb N$, $n\le s_1,s_2\le n+k$, $k:=[\sqrt{n}]$ and $N=2n+k$, the following estimates hold:
$$
\frac{\|P^Ne_1\|}{\|P^Ne_{2s_1}\|}\le e^{-\beta n^2},\ \ e^{-\gamma n^{3/2}}\le\frac{\|P^Ne_{2s_1}\|}{\|P^Ne_{2s_2}\|}\le e^{\gamma n^{3/2}}
\tag3.10
$$
for some positive constants $\beta$ and $\gamma$ and all $n$ large enough (in order to simplify the notations, we will write below $\sqrt{n}$
instead of its integer part).
\endproclaim
\demo{Proof} Indeed, let $\delta_n:=(\lambda_n+\lambda_{n+1})T/2$ and $\delta_0:=T\lambda_1/2$. Then, due to \thetag{2.5} and \thetag{2.6},
$$
\|P^Ne_1\|=e^{-\sum_{l=1}^{2N+1}\delta_l},\ \ \|P^Ne_{2s}\|=e^{-\sum_{l=0}^{2s-1}\delta_l-\sum_{l=1}^{2(N-s)+1}\delta_l}.
\tag3.11
$$
for all $n\le s\le n+k$.
Since
$$
\multline
\sum_{l=0}^{2s-1}\delta_l+\sum_{l=1}^{2(N-s)+1}\delta_l=\sum_{l=1}^{2N+1}\delta_l-\sum_{l=2(N-s)+2}^{2N+1}\delta_l+\sum_{l=0}^{2s-1}\delta_l=\\=
\sum_{l=1}^{2N+1}\delta_l-\sum_{l=0}^{2s-1}[\delta_{l+2(N-s)+2}-\delta_l],
\endmultline
\tag3.12
$$
then, due to \thetag{2.3}, we have
$$
\frac{\|P^Ne_1\|}{\|P^Ne_{2s_1}\|}=e^{-\sum_{l=0}^{2s-1}[\delta_{l+2(N-s)+2}-\delta_l]}\le e^{-\sum_{l=0}^{2n}[\lambda_{2n+l}-\lambda_l]}\le e^{-\beta n^2}
$$
for some positive $\beta$. Thus, the first inequality of \thetag{3.11} is verified. Let us now check the second one. To this end, we transform the left-hand side of \thetag{3.12} as follows:
$$
\sum_{l=0}^{2s-1}\delta_l+\sum_{l=1}^{2(N-s)+1}\delta_l=\sum_{l=0}^{2n}\delta_l+\sum_{l=1}^{2n-1}\delta_l+
\sum_{l=2n+1}^{2s-1}\delta_l+\sum_{l=2n}^{2(N-s)+1}\delta_l
\tag3.13
$$
which implies the second estimate of \thetag{3.11} since two first terms in the right-hand side of \thetag{3.13} are independent of $s$ and the third and fourth sums contain no more than $\sqrt{n}$ terms each of them is not greater than $\delta_{N}\sim Cn$. Lemma 3.3 is proved.
\enddemo
According to Lemma 3.3, for every $n\in\Bbb N$ and $k$ and $N(n)$ as in the statement of the lemma, there exist numbers $A_s(n)$, $0\le s\le k$, such that
$$
1\ge A_s(n)\ge e^{-2\gamma n^{3/2}}
\tag3.14
$$
and
$$
\|A_0(n)P^Ne_{2n}\|=\|A_1(n)P^N e_{2(n+1)}\|=\cdots=\|A_k(n)P^Ne_{2(n+k)}\|= B(n)
\tag3.15
$$
where
$$
e^{-\gamma_2 n^2}\le B(n)\le e^{-\gamma_1 n^2}
\tag3.16
$$
for some positive $\gamma_i$.
\par
We are now ready to describe the perturbation of \thetag{2.23} which will produce the infinite log-doubling factor. To this end, we recall that
the 3D perturbation $\Cal R$ vanishes if $w_1<2$. By this reason, without loss of generality, we may assume that $\Cal R(v_s(t))\equiv 0$ for all $s\in[1-\kappa,1]$, $-\kappa\le t\le0$.
\par
We now split the interval $s\in[1-\kappa,\kappa]$ on infinitely many pieces $I_n$ such that the length $|I_n|=\kappa 2^{-n}$. Further, every segment
$I_n$ will be further divided on $2^{\sqrt{n}}$ equal subintervals $I_{n,p}$, $p=1,\cdots, \sqrt{n}$. Then, of course, $|I_{n,p}|\ge \kappa2^{-2n}$.
For every interval $I_{n,p}$, we now fix $s_{n,p}\in I_{n,p}$, for instance, the midpoint of that interval the corresponding trajectory $v_{s,p}(t)$
in a such way that
$$
v_{n,p}(0)=s_{n,p}\in I_{n,p}.
$$
At the next step, we introduce a family of smooth cut-off functions $\psi_{n,p}(x,y,w_1)$ such that
$$
\psi_{n,p}(v_{n,p}(t))\equiv 1,\ t\in[-\kappa,0],\ \ \psi_{n,p}(v_{n_1,p_1}(t))\equiv0,\ (n_1,p_1)\ne(n,p)
\tag3.17
$$
Obviously, such functions exist. Moreover, since $|I_{n,p}|\ge \kappa2^{-2n}$, we may fix these function in a such way that
$$
\|\psi_{n,p}\|_{C^R(\Bbb R^3)}\le M_R 2^{2Rn},\ \ \text{for all}\ R\in\Bbb N,
\tag3.18
$$
where the constants $M_R$ are independent of $p$ and $n$.
\par
We also need one more smooth non-zero bump  function $\theta(x)\in C_0^\infty(\R)$ such that
$$
\theta(x)\ge0,\ \ \theta(x(t))\equiv0 \ \ \text{for all}\ t\in[-T,0],\ t\notin(-\kappa,0)
\tag3.19
$$
and a family of constants
$$
K_{n}:=\int_{-\kappa}^0e^{-\lambda_{2n}h}\theta(x(h))\,dh.
\tag3.20
$$
Clearly,
$$
1\ge K_n\ge 2^{-C n}
\tag3.21
$$
for the properly chosen positive constant $C$.
\par
Let
$$
g_p(n)\in\sppan\{e_{2(n+1)},\cdots,e_{2(n+\sqrt{n})}\},\ p=1,\cdots, 2^{\sqrt{n}}
$$
be the vortices of the $\sqrt{n}$-dimensional unit cube $[0,1]^{\sqrt{n}}$ enumerated in some order. Finally, we want to construct a smooth perturbation $T_n(x,y,w_1)\in \Cal L(H,H)$ in such way that
$$
\multline
v_{n,p}(2(2n+\sqrt{n})T)=(x(0),y(0),w_{n,p}(2(2n+\sqrt{n})T),\\
w_{n,p}(2(2n+\sqrt{n})T)=\\=(w_{n,p}(2(2n+\sqrt{n})T),e_1)e_1+e^{-\beta n^2/2}B(n)g_p(n),\ \ n\in\Bbb N,\ \ p=1,\cdots,\sqrt{n}
\endmultline
\tag3.22
$$
This can be done by the following formula:
$$
\multline
T_n(x,y,w_1) w=\\= 2^{-\beta n^2/2}\sum_{k=1}^{\sqrt{n}}\sum_{p=1}^{2^{\sqrt{n}}}\theta(x)\psi_{n,p}(x,y,w_1)K_{n+k}^{-1}A_k(n)g_p(n)_k w_{2(n+k)}e_{2(n+k)},
\endmultline
\tag3.23
$$
where $g_p(n)_k\in\{0,1\}$ denotes the $e_{2(n+k)}$ coordinate of the vortex $g_p(n)$. Indeed, \thetag{3.22} follows directly
 from the definition of the constants $A_k(n)$, $B(n)$, $K_m$, the cut-off functions and the fact that, for the $e_{2(n+k)}$ coordinate of the vortex
 $g_p(n)$, we factually need to solve the equation
 $$
 \frac d{dt}w_{2(n+k)}+\lambda_{2(n+k)}w_{2(n+k)}=e^{-\beta n^2/2}K_{2(n+k)}^{-1}A_k(n)g_p(n)_k\theta(x(t)),
 $$
 with $w_{2(n+k)}(-\kappa)=0$ (we recall that $\psi_{n,p}(v_{n,p}(t))\equiv1$). Moreover, due to \thetag{3.23}, \thetag{3.21} and \thetag{3.18},
 we conclude that
 $$
 \|T_n\|_{C^R(\R^3,\Cal L(H,H))}\le M'_Re^{-\beta n^2/2}2^{kRn},
 \tag3.24
 $$
 where the positive constant $k$ is independent of  $n,R\in\Bbb N$ and $M'_R$ depends only on $R$. Finally, we define
 $$
 T(x,y,w_1)w:=\sum_{n=n_0}^\infty T_n(x,y,w_1)w,
 \tag3.25
$$
where $n_0$ is a sufficiently large number. Then, using \thetag{3.24}, we see that the perturbation $T$ is $C^\infty$-smooth and for any $R\in\Bbb N$,  the $C^R$-norm of it  can be made arbitrarily small by fixing $n_0$ being large enough.
\par
We now claim that the attractor $\Cal A$ of the perturbed version
$$
\multline
\frac d{dt} x=f(x,y)+R_1(x,y,w_1),\ \frac d{dt} y(t)=g(x,y)+R_2(x,y,w_1),\\  \frac d{dt} w+Aw=\Cal R(x)w+R_3(x,y,w_1)e_1+T(x,y,w_1)w
\endmultline
\tag3.26
$$
of \thetag{2.25} has infinite Log-doubling factor. The proof of this fact is based on \thetag{3.22} and \thetag{3.10}. Indeed, according to \thetag{3.22}, the attractor $\Cal A$ contains a sequence of vortices of "almost cubes"
$$
u(0)+(w_{n,p}(2(2n+\sqrt{n})T),e_1)e_1+\eb_n g_p(n)\in\Cal A, \ \ n\ge n_0, \ p=1,\cdots 2^{\sqrt{n}}
\tag3.27
$$
with $\eb_n:=e^{-\beta n^2/2}B(n)$ and $u(0)=(x(0),y(0),0)$. Furthermore, due to the first estimate of \thetag{3.10}, estimates \thetag{3.14} and \thetag{3.21} and the definition \thetag{3.23}, we see that
$$
|(w_{n,p}(2(2n+\sqrt{n})T),e_1)|\le e^{-\beta n^2/2}e^{Cn^{3/2}}\eb_n\ll\eb_n
\tag3.28
$$
if $n_0$ is large enough. Thus, the impact of the second term in \thetag{3.27} is neglecting for what follows and, with the minor loss of rigoricity, we may assume that
$$
u(0)+\eb_n g_p(n)\in\Cal A
\tag3.29
$$
for all $n\ge n_0$ and all $p=1,\cdots,2^{\sqrt{n}}$. As the following lemma shows, that is sufficient to conclude that the log-doubling dimension of $\Cal A$ is infinite.
\proclaim{Lemma 3.4} Let $\Cal B\subset H$ be a compact set in our Hilbert space $H$ such that
$$
\eb_n g_p(n)\in\Cal B
\tag3.30
$$
for all $n\ge n_0$ and all $p=1,\cdots,2^{\sqrt{n}}$ and $\eb_n\sim e^{-\gamma n^2}$ (with some positive $\gamma$). Then the
log-doubling factor of $\Cal B$ is infinite.
\endproclaim
\demo{Proof} Indeed, consider the ball $B_{r_n\eb_n}$ of radius $r_n:=n^{1/4}/2$ centered at the center of the cube spanned by $\eb_n g_p(n)$, $p=1,\cdots,2^{\sqrt{n}}$. Then, obviously,
$$
\eb_n g_p(n)\in B_{r_n\eb_n}\cap\Cal B,
\tag3.31
$$
but
$$
N_{\eb_n/2}(B_{r_n\eb_n}\cap \Cal B,\Cal B)\ge 2^{\sqrt{n}}.
\tag3.32
$$
Using now the obvious estimate
$$
N_{\eb_n/2}(B_{r_n\eb_n}\cap \Cal B,\Cal B)\le D_{\eb_n/2}(\Cal B)^{\log_2 r_n},
\tag3.33
$$
we see that
$$
\log_2 D_{\eb_n/2}(\Cal B)\ge\sqrt{n}/\log_2 r_n\ge 1/2\sqrt{n}\sim \frac12\(\frac1\gamma\log\frac2{\eb_n}\)^{1/4}.
\tag3.34
$$
Thus, the log-doubling dimension of $\Cal B$ is indeed infinite and the lemma is proved.
\enddemo
The fact that the global attractor $\Cal A$ of system \thetag{3.26} follows immediately from \thetag{3.29} and the proved lemma. Theorem 3.2 is proved.
\enddemo

\head Section 4. "Bad" projections and the case of finite smoothness
\endhead
In this section, we construct several related counterexample showing that the fractal dimension of the global attractor may depend on the choice of the phase space in the case where the nonlinear operator $F(u)$ is not $C^\infty$-smooth as well as the example of the $C^\infty$-smooth attractors whose image under the orthogonal projector of the finite rank has infinite log-doubling factor and, by this reason, cannot be embedded into any finite-dimensional log-Lipschitz manifold.
\par
Both examples will be based on the same construction which we described below. Namely, let us consider the following 2D system:
$$
\cases
x'=-x(x^2+y^2-1),\\
y'=-y(x^2+y^2-1).
\endcases
\tag4.1
$$
This system can be rewritten in polar coordinates $x+i y=Re^{i\varphi}$ as follows:
$$
R'=-R(R^2-1),\ \ \varphi'=0
\tag4.2
$$
Thus, the global attractor $\Cal A_0\subset\R^2$ of this system consists of all points $R\in[-1,1]$, $\varphi\in[0,2\pi]$ and is $C^\infty$-diffeomorphic to a cone.
\par
We now split the segment $\varphi\in[0,2\pi]$ on the infinite number of intervals $I_n$, of lengthes $|I_n|=A_n$, $n=1,2\cdots$ satisfying
$$
\sum_{n=1}^\infty A_n<2\pi,
\tag4.3
$$
fix $\varphi_n\in I_n$ and a family of the cut-off functions $\psi_n(\varphi)\in C_0^\infty(\R)$ such that $\psi_n(\varphi_k)=\delta_{nk}$ satisfying
$$
\|\psi_n\|_{C^k(\R)}\le C_k A_n^{-k},
\tag4.4
$$
where the constant $C_k$ is independent of $n$, and the cut-off function $\theta(R)\in C_0^\infty(\R)$ such that
$$
\theta(R)\equiv0,\ \ R\le0,\ \ \theta(R)\equiv1,\ \ R\in[1/2,1].
\tag4.5
$$
Finally, we introduce the operator $F: \R^2\to H$ via
$$
F(x,y):=\sum_{n=1}^\infty B_n\theta(x,y)\psi_n(x,y)e_n,
\tag4.6
$$
where the monotone decreasing to zero sequence $B_n$ will be specified below, and consider the coupled
system
$$
\cases
x'=-x(x^2+y^2-1),\ \ y'=-y(x^2+y^2-1),\\
\Dt v+A v=F(x,y)
\endcases
\tag4.7
$$
The elementary properties of this system are collected in the following lemma:

\proclaim{Lemma 4.1} Let $H^s:=D(A^{s/2})$ be the scale of H-spaces generated by the operator $A$. Then,
\par
1) The map $F$ defined by \thetag{4.6} belongs to
$F\in C^k(\R^2,H^s)$ iff
$$
\sup_{n\in\Bbb N}\{B_n\lambda_n^s A_n^{-k}\}<\infty.
\tag4.8
$$
\par
2) Equation \thetag{4.7} possesses a global attractor $\Cal A$ in $H$ which contains the following sequence of points:
$$
P_n:=(\cos\varphi_n,\sin\varphi_n,B_n\lambda_n^{-1}e_n)\in\Cal A
\tag4.9
$$
for $n=1,2,\cdots$. In addition, let $Q_2:\R^2\times H\to H$ be the orthoprojector to the $v$-component of \thetag{4.7}. Then, the projection
$Q_2\Cal A$ contains the following sequence of segments:
$$
S_n:=\{se_n,\ \ s\in[0,B_n\lambda_n^{-1}]\}\subset Q_2\Cal A.
\tag4.10
$$
\endproclaim
\demo{Proof} Indeed, the first assertion of the lemma is an immediate corollary of the definition \thetag{4.6}, estimate \thetag{4.4} and the fact that, for every fixed $(x,y)\in\Bbb R^2$, at most one term in \thetag{4.6} is non-zero. To verify \thetag{4.9} it remains to note that, due to the definition of the cut-off functions and the operator $F$, $P_n$ is an equilibrium of \thetag{4.7} for every $n\in\Bbb N$ and \thetag{4.10} follows from the fact that there is a heteroclinic orbit connecting $(\cos\varphi_n,\sin\varphi_n,0)$ and $P_n$, again for every $n\in\Bbb N$.
\enddemo
Note that the embedding \thetag{4.10} already shows that
$$
\lim_{\eb\to0}D_\eb(Q_2\Cal A)=\infty
\tag4.11
$$
if the decreasing sequence $B_n$ is strictly positive and, by this reason, $Q_2\Cal A$ cannot be embedded into any finite-dimensional Lipschitz manifold, see \cite{24}. The next theorem shows that, under the proper choice of the sequence $B_n$, the embedding into the Log-Lipschitz manifolds is also impossible.

\proclaim{Theorem 4.2} Let $\lambda_n\sim n^\kappa$ (which corresponds to the case where $A$ is an elliptic  differential operator in a bounded domain) and let $B_n\sim e^{-\sqrt{n}}$. Then the non-linearity $F$ in \thetag{4.7} belongs to $C^\infty(\R^2,H^s)$ for any $s\in\R$, but
$$
\dim_{Log-D}(Q_2\Cal A)=\infty
\tag4.12
$$
and, by this reason, $Q_2\Cal A$ cannot be embedded into any Log-Lipschitz manifold.
\endproclaim
\demo{Proof} Indeed, the smoothness of $F$ follows from the criterium \thetag{4.8}. Let us check \thetag{4.12}. According to \thetag{4.10},
for any $n\in\Bbb N$, the points $B_n\lambda_n^{-1}e_k$, $k=1,\cdots,n$ belong to $Q_2\Cal A$ and, therefore,
$$
D_{\eb_n}(Q_2\Cal A)\ge n,\ \ \eb_n=\lambda_n^{-1}B_n
\tag4.13
$$
Since $\eb_n\sim n^{-\kappa}e^{-\sqrt{n}}$, we see that $n\sim(\log\frac1{\eb_n})^2$ and \thetag{4.12} holds.
\enddemo
\remark{Remark 4.3} Note that the equality \thetag{4.12} does not hold for
 the attractor $\Cal A$ itself, but only for its "bad" projection $Q_2\Cal A$. In addition, modifying equation \thetag{4.7} as follows
 $$
\cases
x'=-\beta x(x^2+y^2-1),\ \ y'=-\beta y(x^2+y^2-1),\\
\Dt v+A v=\beta F(x,y)
\endcases
\tag4.14
$$
where $\beta\ll1$ is a small positive parameter, we see that the assertion of the Theorem 4.2 still holds, but now the spectral gap condition is satisfied and there is a $C^1$-smooth inertial manifold diffeomorphic to $\R^2\ni(x,y)$ containing the global attractor $\Cal A$ (by decreasing $\beta>0$ this manifold can be made of the class $C^k$ for any fixed $k>0$).
\endremark
As the last example, we now consider the case where the operator $F$ has only finite smoothness. We first remind that, in the case of infinite smoothness, say, $F\in C^\infty(\R^2,H^s)$, for any $s\in\R_+$, due to the parabolic  smoothing property, the solution operator $S(t):\R^2\to H\to\R^2\to H$ of equation \thetag{4.7} satisfies
$$
\|S(1)u_1-S(1)u_2\|_{\R^2\times H^s}\le C_s\|u_1-u_2\|_{\R^2\times H},\ \ u_i=(x_i,y_i,v_i)\in\Cal A,\ \ i=1,2
\tag4.15
$$
and, by this reason, the fractal dimension of the attractor $\Cal A$ is {\it the same} in all spaces $H^s$. However, in the case where $F$ is only finite smooth, there is a limit exponent $s_0$ (the largest $s$ for which \thetag{4.15} holds) and beyond this exponent we cannot guarantee that the fractal dimension remains the same. As the next theorem show, the dimension may indeed start to grow and even become infinite for $s>s_0$.
\proclaim{Theorem 4.4} Let
$$
\lambda_n\sim n^2, \ \ A_n\sim \frac{1}{\lambda_n^{1/2}(\log\lambda_n)^2},\ \ B_n=\frac{A_n}{\log\lambda_n}.
\tag4.16
$$
Then, $F\in C(\R^2,H^1)\cap C^1(\R^2,H)$ and the fractal dimension $\dim_f(\Cal A,\R^2\times H^s)$  of the attractor $\Cal A$ of \thetag{4.14} (with $\beta\ll1$) equals to $2$ for $s\le2$ and starts to grow when $s>2$. Moreover, $\Cal A$ is still compact in $\R^2\times H^3$, but its fractal dimension is infinite there:
$$
\dim_f(\Cal A,\R^2\times H^3)=\infty.
\tag4.17
$$
\endproclaim
\demo{Proof} Indeed, the regularity of the map $F$ follows from the criterium \thetag{4.8}, the compactness of the attractor in $H^3$ is a straightforward corollary of the fact that $\lambda_n^{3/2}\lambda_n^{-1}B_n\to0$ as $n\to\infty$ and the explicit construction of the attractor. The fact that the fractal dimension $\dim_f(\Cal A,\R^2\times H)=2$ follows from the existence of a 2D inertial manifold for that problem (spectral gap assumption is satisfied if $\beta$ is small enough) plus the existence of the 2D unstable manifold belonging to the attractor. The fact that this dimension remains equal to $2$ in the spaces $\R^2\times H^s$ with $s\le2$ follows from the validity of the smoothing property \thetag{4.15} for $s\le2$ (the limit exponent $s_0=2$ here). Thus, it only remains to verify \thetag{4.17}.
\par
Since the fractal dimension does not increase under the orthogonal projection, \thetag{4.10} gives
$$
\multline
\dim_f(\Cal A,\R^2\times H^3)\ge\dim_f(Q_2\Cal A,H^3)\ge \dim_f(\cup_{n=1}^\infty S_n,H^3)\ge\\\ge \dim_f(\cup_{n=1}^\infty\{\lambda_n^{-1}B_ne_n\},H^3)=\dim_f(\cup_{n=1}^\infty\{\lambda_n^{-3/2}\frac1{(\log \lambda_n)^3}e_n\},H^3)=\\
\dim_f(\cup_{n=1}^\infty\{\frac1{(\log n)^3}e_n\},H)=\infty
\endmultline
\tag4.18
$$
and the theorem is proved.
\enddemo
\remark{Remark 4.5} Using the fact that the countable union of the sets with Hausdorff dimension $2$ also has the Hausdorff dimension $2$,
one can verify that the Hausdorff dimension of the attractor equals $2$ in all $H^s$, $s\le3$
$$
\dim_H(\Cal A, \R^2\times H^3)=2,\ \ s\le3
\tag4.19
$$
and, in particular, for $s=3$, we have finite Hausdorff and infinite fractal dimension.
\endremark

\head Section 5. Concluding remarks and open problems
\endhead
In this concluding section, we give some remark clarifying the counterexamples given above and indicate some related open problems.

\remark{Remark 5.1} Note that, most of our counterexamples are constructed under the assumption that the nonlinear function $F$ is $C^\infty$-smooth.
It would be interesting to give the corresponding examples in the class of real {\it analytic} functions. In particular, it is not difficult to see that the example of Theorem 1.1 which gives the absence of a $C^1$-inertial manifold can be straightforwardly realized in the class of analytic $F$. However, extending the examples of Sections 2 and 3 to analytic nonlinearities is not so straightforward and could be an interesting problem.
\par
Another question is how generic is the situation when the $C^1$, Lipschitz or Log-Lipschitz inertial manifold does not exist. Again, the $C^1$-situation is clear (indeed, the non-existence result of Theorem 1.1 is based on the robust spectral properties of two equilibria which will survive under small perturbations and this shows that for any close equation of the form \thetag{0.4}, we will have the non-existence result) and it is again not clear whether or not the same is true for the Lipschitz and Log-Lipschitz cases.
\endremark
\remark{Remark 5.2} We need to emphasize that the non-smoothness of the homeomorphism $P:\Cal A\to P\Cal A$ {\it does not} imply in general that the conjugated DS $\widetilde S(t)=P S(t)P^{-1}$ is also not smooth. Indeed, there are many examples in the modern DS theory when two smooth DS are conjugated by the non-smooth homeomorphism only (see e.g., Grobman-Hartman type theorems, \cite{11,14}). Thus, the proved non-existence results for the Lipschitz and Log-Lipschitz Mane projections do not imply automatically that the dynamics on the attractor is "infinite-dimensional" or "non-smooth", but only indicates that potentially smooth and finite-dimensional dynamics is {\it embedded} into the infinite-dimensional space in the essentially non-smooth way. Moreover, in all examples of this paper the factual dynamics on the attractor is {\it trivial} and at least intuitively smooth and low dimensional (e.g., one stable limit circle plus several hyperbolic equilibria with instability indexes not larger than 2).
\par
On the other hand, the structure of the Poincare map of this limit circle is similar to the standard shift operator on the space $l_2(\Bbb Z)$, see Proposition 2.1 and it is not easy to interpret  this as a "finite-dimensional" phenomenon. Thus, it would be interesting to see the impact of such and similar structures not only to the regularity of the attractor's embedding, but also to the dynamical properties, bifurcations, etc. For instance, it would be important to construct the examples showing that the dynamics on the attractor cannot be conjugated even topologically to the {\it smooth} dynamics on compact set of $\R^N$ (or to prove that it is always possible). Note that, due to the superexponential convergence to the limit circle, see Section 2 and 3, we see that this conjugating homeomorphism (if exists) cannot be  H\"older continuous.
\endremark
\remark{Remark 5.3} The constructed counterexamples, in a sense, show the {\it limitations} of the abstract functional model approach and indicate that many results on the abstract equation \thetag{0.4} obtained so far are sharp and there is no hope to improve them essentially. For instance, we know that the Mane projections of the attractor are always bi-H\"older  with the H\"older exponent arbitrarily close to one, but as we have seen, not bi-Lipschitz or bi-Log-Lipschitz in general. We also have seen that the known obstacles to the existence of the $C^1$-inertial manifolds (see Section 1) can also forbid the existence of Lipschitz and Log-Lipschitz manifolds and, at least on the level of abstract equations, the spectral gap condition is responsible also for that Lipschitz and Log-Lipschitz manifolds. In addition, from the log-convexity and backward uniqueness theorems, see \cite{1, 15}, we know that the following $1/2$-Log-Lipschitz continuity holds on the attractor
$$
\|Au_1-Au_2\|_{H}\le C\|u_1-u_2\|_H\(\log\frac C{\|u_1-u_2\|_H}\)^{1/2},\ \ u_1,u_2\in\Cal A
\tag5.1
$$
and the example from Section 2 shows that we cannot expect better exponent than $1/2$ in this relation (exactly the exponent $1/2$ gives the superexponential growth/decay at the form of $e^{-t^2}$).
\par
However, the obtained counterexamples {\it do not imply} at least in a straightforward way that the idea with Lipschitz and Log-Lipschitz Mane projections will not work for the {\it concrete} classes of equations of mathematical physics, like reaction-diffusion systems,
or 2D Navier-Stokes equations on a torus. Indeed, in order to extend the above counterexamples from the abstract problem \thetag{0.4} to the concrete classes of equations, say, the reaction-diffusion ones, as the most difficult (from our point of view, of course) step, we need to build up the equilibria $u_1$ and $u_2$ of the equation
$$
\Dt u=a\Delta_x u-f(x,u)
\tag5.2
$$
with the spectral properties described in Section 1 and this leads to the  non-trivial and not properly understood {\it inverse spectral} problem.
We also mention that, in contrast to the abstract situation, the spectral gap condition  may be not necessary. Indeed, for the 1D equation \thetag{0.7} mentioned in Introduction, the spectral gap condition is violated (due to the presence of $u_x$ in the non-linearity). However, the linearized operator
$$
L_u v:=v_xx-f'_u v-f'_{u_x}v_x
$$
is similar to the self-adjoint one $w_{xx}-K(x)w$ (with $w=e^{1/2\int f'_{u_x}\,dx}v$) for which the spectral gap condition holds and, therefore, the counterexample of Section 1 is impossible. This property has been used in a crucial way in the Romanov's proof of the existence of bi-Lipschitz Mane projections for that case.
\par
Thus, it would be interesting to find other obstacles which make the spectral construction of Section 1 impossible.
\endremark

\enddocument